\documentclass{elsarticle}
\usepackage{hyperref}
\usepackage{mathrsfs}
\let\memoldbibsection\bibsection
\let\bibsection\relax
\usepackage{amsrefs}         
\let\bibsection\memoldbibsection
\usepackage{amsmath}
\def\cl{\mathop{\rm cl}\nolimits}

\def\int{\mathop{\rm int}\nolimits}

\def\RO{\mathop{\mathsf{RO}}\nolimits}
\def\R{\mathop{\mathsf{R}\hspace{0mm}}\nolimits}

\makeatletter
\def\ps@pprintTitle{%
  \let\@oddhead\@empty
  \let\@evenhead\@empty
  \let\@oddfoot\@empty
  \let\@evenfoot\@oddfoot
}
\makeatother

\usepackage{amsfonts}
\usepackage{amssymb}
\newtheorem{thm}{Theorem}
\newtheorem{lem}[thm]{Lemma}
\newtheorem{exa}{Example}
\newtheorem{pro}{Proposition}
\newdefinition{rmk}{Remark}
\newtheorem{cor}{Corollary}
\newtheorem{question}{Question}
\newtheorem{definition}{Definition}
\newproof{pf}{Proof}
\begin{document}
\begin{frontmatter}
\title{On selectively highly divergent spaces} 
\author[1]{Carlos David Jiménez-Flores\fnref{fn1,fn5}}
\ead{carlosjf@ciencias.unam.mx}
\author[1]{Alejandro Ríos-Herrejón\fnref{fn2}}
\ead{chanchito@ciencias.unam.mx}
\author[2]{Alejandro Darío Rojas-Sánchez}
\ead{adrojas@up.edu.mx}
\author[1]{Elmer Enrique Tovar-Acosta\fnref{fn4}}
\ead{elmer@ciencias.unam.mx}
\affiliation[1]{organization={Universidad Nacional Autónoma de México.},
            addressline={Av. Universidad 3000}, 
            city={CDMX},
            postcode={4510}, 
            country={México}}  
\affiliation[2]{organization={Universidad Panamericana},
            addressline={Augusto Rodin No. 498, Col. Insurgentes Mixcoac}, 
            city={CDMX},
            postcode={03920}, 
            country={México}} 
\fntext[fn1]{The first author acknowledges financial support from CONACyT, grant no. 816193}
\fntext[fn2]{The second author acknowledges financial support from CONACyT, grant no. 814282.}
\fntext[fn4]{The fourth author acknowledges financial support from CONACyT, grant no. 829699}
\fntext[fn5]{This paper is a part of the doctoral research of the first author.}
\begin{abstract}
We say that a topological space $X$ is selectively highly divergent (SHD) if for every sequence of non-empty open sets $\{U_n: n\in\omega \}$ of $X$, we can find $x_n\in U_n$ such that the sequence $(x_n)$ has no convergent subsequences. We investigate the basic topological properties of SHD spaces and we will exhibit that this class of spaces is full of variety. We present an example of a SHD space wich has a non trivial convergent sequence and with a dense set with no convergent sequences. Also, we prove that if $X$ is a regular space such that for all $x\in X$ holds $\psi(x,X)>\omega$, then $X_\delta$ (the $G_\delta$ modification of $X$) is a SHD space and, moreover, if $X$ homogeneous, then $X_\delta$ is also homogeneous. Finally, given $X$ a Hausdorff space without isolated points, we construct a new space denoted by $sX$ such that $sX$ is extremally disconnected, zero-dimensional Hausdorff space, SHD with $|X|=|sX|$, $\pi w(X)=\pi w(sX)$ and $c(X)=c(sX)$ where $\pi w$ and $c$ are the cardinal functions $\pi$-weight and celullarity respectively. 
\end{abstract}
\begin{keyword}
\MSC[2010] 54A20 \sep 54A25 \sep 54B20 \sep 54D40 \sep 54G05\\
sequences \sep compactifications \sep remainder \sep realcompact \sep absolutes \sep $F$-spaces \sep isolated points
\end{keyword}
\end{frontmatter}
\section{Introduction}
The properties of sequences have always been a subject of study, although they do not generally completely characterize the topology of a space. A couple of relevant properties, as highly divergente sequences, are those appearing in \cite{altamentedivergente}. On the other hand, A. Dorantes-Aldama and D. Shakhmatov define a topological property for sequences as it appears in \cite[Definition~2.1]{Dorantes}, and based on this, they stated the term selectively $S$ space in \cite[Definition~2.3(iii)]{Dorantes} with $S$ a topological property for sequences. The property of being selectively highly divergent arises from considering the ``highly divergent" property together with the definition of a selectively $S$ space. It should be noted that the highly divergent property does not satisfy \cite[Definition~2.1]{Dorantes}. As a result, a very diverse class of spaces with quite relevant properties was obtained, which will be the object of study throughout this article.\\\\All topological spaces are assumed to have no separation unless otherwise stated. Additionally, the notation and definitions we use can be found in \cite{Porter}, \cite{Gillman}, and \cite{Engelking}. The notation and definition for cardinal functions is as in \cite[Chapter~1]{handbook}. Here, as usual, $X^{*}$ denotes the space $\beta X\setminus X$, i.e., the remainder of the Stone-Čech compactification of $X$.
\section{Basic properties}
Inspired by \cite[Definition~2.3]{Dorantes} and \cite[Definition~2]{altamentedivergente}, we say that a topological space $X$ is \textit{selectively highly divergent} (SHD from here for short) if for every sequence of non-empty open sets $\{U_n: n\in\omega \}$ of $X$, we can find $x_n\in U_n$ such that the sequence $(x_n)$ has no convergent subsequences. Clearly, this is a topological property that is inherited by open subspaces. Furthermore, it follows from the definition that if a topological space $X$ has a point of countable character, then it cannot be SHD and thus, no metrizable space is SHD. The following remark will be relevant for the subsequent development of this work:
\begin{rmk}
Let $(X,\tau)$ be a topological space and $\mathcal{B}$ a $\pi$-base for $X$. If for any sequence of open sets $\{U_n : n \in \omega\}\subseteq \mathcal{B}$, it is true that for each $n \in \omega$ there exists $x_n \in U_n$ such that the sequence $(x_n)$ has no convergent subsequences, then $X$ is SHD.
\end{rmk}

The following theorem is easy to prove and establishes the relationship between the property of being SHD and topological products:
\begin{thm} \label{productomamado}
If $\{X_i : i \in I\}$ is a family of topological spaces and at least one of them is SHD, then $\prod\{X_i : i \in I\}$ is SHD.
\end{thm}

We obtain an analogous result for the topological sum:
\begin{thm}
If $\{X_i: i\in J\}$ is a family of topological spaces, and all of them are SHD, then $\bigoplus \{X_i: i\in J\}$ is SHD.
\end{thm}
\begin{pf}
Let $X=\bigoplus \{X_i: i\in J\}$. First note that if a sequence $(z_n)$ converges in $X$, then there exist some $j\in J$ such that a tail of the sequence is contained in $X_j$, i.e., convergente sequences in $X$ are eventually cointained in one and only one of the spaces. Let $\{U_n: n\in\omega \}$ be a family consisting of basic open sets, i.e., every $U_n$ is an open set of some $X_i$. Let $\alpha$ the function that takes a natural number $n\in\omega$ and asigns an element $j\in J$ with the property that $U_n\subseteq X_j$. Thanks to the fact that the spaces of the topological sum are disjoint by pairs, the index obtained by $\alpha$ is unique. Then, for every $n\in\omega$, we'll write $\alpha(n)$ instead the index $j\in J$ and, from the previous argument, $\alpha(n)\in J$ is the only element in $J$ with the property that $U_n\subseteq X_{\alpha(n)}$. Consider the following cases:\\\\
		\textbf{Case 1.} For every $j\in J$, the set $\alpha^{-1}[j]$ is finite. In this case, simply select any point $x_n\in U_n$. Note that $(x_n)$ doesn't have convergent sequences, since every subsequence is infinitely oscillating between the  different $X_{\alpha(n)}$, and so it can't converge by the observation made at the start of the proof.\\
		\textbf{Case 2.} There exists $j\in J$ such that $\alpha^{-1}[j]$ is infinite. Let $A=\{ j \in J : |\alpha^{-1}[j]|=\aleph_0 \}$. For each $i\in A$, $\alpha^{-1}[i]$ is a infinite subset of $\omega$, say $\{ n(k,i) : j\in \omega\}=\alpha^{-1}[i]$, where $n(k,i)<n(k+1,i)$ for every $k\in \omega$. In this notation, we have that $\{ U_{n(k,i)} : k\in\omega\}$ is a countable family of open sets in $X_i$, and since $X_i$ is a SHD space, we can take $x_{n(k,i)}\in U_{n(k,i)}$ such that $(x_{n(k,i)})_{k\in\omega}$ is a sequence with no convergente subsequences.\\   Let $B=\omega\setminus \bigcup_{i\in A} \alpha^{-1}[i]$. For every $n\in\omega$ simply take $x_n\in U_n$. Define a new sequence in the following way:
\begin{equation*}
y_n=\left\{
\begin{array}{ccc}
x_{n(k,i)} & \text{if} & n=n(k,i) \ \text{for some} \ k\in\omega, i\in A\\\\
x_n & \text{if} & n\in B
\end{array}
\right.
\end{equation*}
Note that since the family $\mathcal{U}=\{ \alpha^{-1}[j] \ | \ j\in K\}\cup \{B\}$ is made up of disjoint sets and $\bigcup\mathcal{U}=\omega$, this is a well defined sequence. Let's prove that none of the subsequences of $(y_n)$ converges. Take a subsequence $(y_{n_k})$. We have the following cases:\\\\
		\textbf{Case 1.} If $|\{ n_k : k\in \omega\}\cap B|=\aleph_0$, we can think the set $ \{ n_k :  k\in \omega\}\cap B $ as a sequence, say $(n_{k_\ell})_{\ell\in\omega}$. Then $y_{n_{k_\ell}}=x_{n_{k_\ell}}$, so $(y_{n_{k_l}})_{\ell\in\omega}=(x_{n_{k_l}})_{\ell\in\omega}$ is a non convergent subsequence of $(y_{n_k})_{k\in\omega}$, and so $(y_{n_k})_{k\in\omega}$ can't converge.\\
		\textbf{Case 2.} If there exists some $i\in A$ such that $|\{ n_k : k\in \omega\}\cap \alpha^{-1}[i]|=\aleph_0$, then, in the same way as the previous case, we can think the set $\{n_k : k\in \omega\}\cap \alpha^{-1}[i]$ as $(n_{k_\ell})_{\ell\in\omega}$. Then $(y_{n_{k_\ell}})_{\ell\in\omega}$ is a subsequence of both $(y_{n_k})_{k\in\omega}$ and $(x_{n(k,i)})_{k\in\omega}$, and since the last one doesn't have convergent subsequences, we conclude that $(y_{n_{k_\ell}})_{\ell\in\omega}$ does not converge, and thus $(y_{n_k})_{k\in\omega}$ also does not converge.\\
		\textbf{Case 3.} If for all $i\in A$, the set $\{ n_k : k\in \omega\}\cap \alpha^{-1}[i]$ is finite, and $\{ n_k : k\in \omega\}\cap B$ is also finite, then, $(y_{n_k})$ can't converge by the observation made at the start of the proof.\\Therefore, we conclude that every subsequence of $(y_n)$ is divergent, and so $X$ is a SHD space.

\end{pf}
Once we have analyzed whether a topological property is hereditary, productive, or additive, it is natural to consider under which class of functions it is preserved. To show that our property of being SHD is in fact rarely preserved, let us consider the following example:
\begin{exa}
Let $X$ be an infinite, compact and SHD space\footnote{In Proposition~\ref{notaalpie} we will present a wide range of SHD spaces with these properties.} and $Y$ to be a non-SHD space (for example, $[0,1]$ as it is first countable). Thus, we know that $X\times Y$ is a SHD space. However, the projection $\Pi_Y:X\times Y\to Y$ is not only continuous, open, and onto, but also a perfect map by Kuratowski's Theorem. But interestingly, despite being almost a homeomorphism, $Y$ is not SHD. This demonstrates that the property of being SHD is not preserved under continuous, surjective, open, and perfect mappings.
\end{exa}
Clearly if $X$ is SHD and $f:X\to Y$ is bijective, continuous and $f^{-1}$ is sequentially continuous, then $Y$ is SHD. Also notice that the SHD space $X$ constructed in Example \ref{Ejemplo2} and the non-SHD space constructed in Example \ref{Ejemplo3}, here denoted by $Y$, have the property that the identity function condenses $X$ onto $Y$.\\\\The following result is not difficult to prove and will be useful in some of our next constructions:
\begin{pro} \label{lemma_SS1}
Let $X$ be a topological space. If every non-empty open subset of $X$ is infinite and $X$ does not admit non-trivial convergent sequences, then $X$ is SHD.
\end{pro}
From \cite[Fact~3.2]{otrosfespacios}, it follows that if $X$ is an $F'$-space without isolated points, the only convergent sequences in $X$ are eventually constant and the open sets are infinite. By Proposition~\ref{lemma_SS1}, $X$ is SHD. These arguments prove the following statement:
\begin{pro} \label{Fespacios}
Let $X$ be a Tychonoff topological space. If $X$ is an $F'$-space without isolated points then $X$ is SHD.
\end{pro}
Since any $F$-space is an $F'$-space, any $F$-space without isolated points is SHD. It is worth noticing that not all SHD spaces are $F$-spaces. To illustrate this, we can refer to a result from \cite[14Q2]{Gillman}, which states that if $X$ and $Y$ are two infinite pseudocompact spaces, then $X\times Y$ is not an $F$-space. Therefore, we can take $Z$ to be a pseudocompact, infinite, SHD space, and consider $Z\times Z$. It follows that $Z\times Z$ is SHD, but it is not an $F$-space.

\section{Constructions and examples}
Considering the definition of a SHD space, one might intuitively assume that a space of this class does not admit non-trivially convergent sequences. However, surprisingly, it is possible to construct a SHD space with non-trivially convergent sequences, along with a couple of other characteristics. An easy way to construct and example is using Theorem \ref{productomamado} and the space $\omega^{*}\times [0,1]$. But the next example has more interesting properties as $X$ is a SHD space which has non-trivial convergent sequences and also has a dense open set in which the only convergent sequences are the trivial ones.
\begin{exa} \label{Ejemplo2}
We know that $Y=\omega^{*}$ is a SHD space since it is an $F$-space with no isolated points. Let $\mathcal{U}=\lbrace U_n \ | \ n\in \omega\rbrace$ be a cellular family consisting of clopen sets in $Y$. Note that $\bigcup_{n\in\omega} U_n$ is SHD as it is an open subset of $Y$. For each $n\in \omega$, select a point $z_n\in U_n$. Furthermore, for each $x\in Y$, fix a local base for that point, denoted by $\mathcal{V}_x=\lbrace V_i \ | \ i \in I \rbrace$, such that $\bigcup \mathcal{V}_x\subseteq U_n$, where $I$ is a sufficiently large set and $x\in U_n$. To simplify notation, we will write $\mathcal{V}_{z_n}$ as $\mathcal{V}_n$. Let $p\notin \beta\omega$. Our space of interest will be $X=\bigcup_{n\in\omega} U_n \cup \lbrace p \rbrace$, endowed with the following topology.\\ 
Let us first consider $\mathcal{B}_n=\lbrace \bigcup_{m\geq n} V_m \ | V_m \in \mathcal{V}_m \rbrace$. Now, let $\mathcal{B}=\bigcup_{n\in\omega} \mathcal{B}_n$. We define the topology by specifying a neighborhood base at each point in the following way:

\begin{itemize}
\item[1)] If $x\in \bigcup_{n\in\omega} U_n$, then the set $\mathcal{V}_x$ remains a neighborhood base for $x$.
\item[2)] A neighborhood base for $p$ is given by the family $\lbrace B\cup \lbrace p\rbrace \ | \ B\in \mathcal{B} \rbrace$.
\end{itemize}

This defines a Hausdorff topology on $X$. Furthermore, this topology is Lindelöf. Now let's examine a couple of important properties of $X$:
\begin{itemize}
	\item[a)] The set $\bigcup_{n\in\omega} U_n$ inherits its original topology. Moreover, $U_n$ remains a clopen set in $X$.
	\item[b)] The sequence $(z_n)$ converges to $p$ in $X$. This follows from the structure of neighborhoods of $p$. From this, we conclude that $Z=\lbrace z_n \ | n \in \omega\rbrace \cup \lbrace p\rbrace$ is a closed set in $X$. Also, note that this set cannot be open, nor is it a neighborhood of $p$.
	\item[c)] Given a non-trivial sequence $(x_n)\subseteq \bigcup_{n\in\omega} U_n$, this sequence cannot converge to any $y\in \bigcup_{n\in\omega} U_n$. This follows from property a) and the fact that in $\bigcup_{n\in\omega} U_n$ there are no non-trivial convergent sequences, as it is a subspace of $\omega^{*}$.
	\item[d)] If $(x_n)\subseteq \bigcup_{n\in\omega} U_n$ is such that for every $n\in \omega$, $x_n\notin \lbrace z_n \ | \ n\in \omega\rbrace$, then $(x_n)$ does not converge to $p$ in $X$.\\
For each $n\in \omega$, let $i(n)$ be the unique natural number such that $x_n\in U_{i(n)}$. We have the following cases:\\\\
		\textbf{Case 1.} The set $\lbrace i(n) \ | \ n\in \omega\rbrace$ is finite. Under these conditions, the sequence is contained in a finite number of the $U_n$ sets and is therefore forced to converge to someone within this finite union, which is a closed set in $X$. Therefore, $(x_n)$ does not converge to $p$.\\
		\textbf{Case 2.} There exists $m\in \omega$ such that the set $A=\lbrace n\in \omega \ | \ i(n)=m\rbrace$ is infinite. If $(x_n)$ were to converge to $p$, it would follow that $\{x_n\}_{n\in A}$ would also converge to $p$, which is impossible since this sequence is contained in the closed set $U_m$.\\
		\textbf{Case 3.} The set $B=\lbrace i(n) \ | \ n\in \omega\rbrace$ is infinite, and for each $m\in \omega$, $i^{-1}[m]$ is finite. For $m \notin B$, we define $W_m = U_m$. On the other hand, if $m\in B$, for each $n\in i^{-1}[m]$, we choose $V_n\in \mathcal{V}_m$ such that $x_n\notin V_n$, and finally we define $W_m=\bigcap_{i(n)=m} V_n$, which is still an open set containing $z_m$ as it is the finite intersection of elements from $\mathcal{V}_m$. To conclude, we define $W=\bigcup_{m\in\omega} W_m\cup \lbrace p\rbrace$. $W$ is an open set in $X$ containing $p$, and it is disjoint from the sequence $(x_n)$. Therefore, $(x_n)$ does not converge to $p$.
	\item[e)] $X$ is SHD. Let $\lbrace W_n \ | \ n \in \omega\rbrace$ be a collection of non-empty open sets in $X$. Without loss of generality, we assume that each $W_n$ is a basic neighborhood. We define $G_n=W_n\setminus Z$, and thus the collection $\lbrace G_n \ | \ n \in \omega\rbrace$ consists of non-empty open sets in $\bigcup_{n\in\omega} U_n$. As this space is SHD, we can choose $x_n\in G_n$ such that $(x_n)$ does not have convergent subsequences in $\bigcup_{n\in\omega} U_n$. Combining this with the previous result, we conclude that none of its subsequences can converge to the point $p$.
\end{itemize}
\end{exa}
\begin{exa} \label{Ejemplo3}
In a similar way to the previous example, we consider the set $X = \bigcup_{n\in\omega} U_n \cup \lbrace p \rbrace$, where the $U_n$ are clopen sets in $\omega^{*}$ and for all $x\in \omega^{*}$ consider $\mathcal{V}_x$ a local basis for $x$ such that $\bigcup\mathcal{V}_x\subseteq U_n$ when $x\in U_n$ and $p \notin \beta \omega$. We define the topology by specifying a neighbourhood base at each point in the following way:
\begin{itemize}
\item For each $x \in \bigcup_{n\in\omega} U_n$, the set $\mathcal{V}_x$ remains a neighborhood base for $x$.
\item For $p$, a neighborhood base is given by the family $\lbrace \bigcup_{m\geq n} U_m \cup \lbrace p\rbrace \ | \ n \in \omega\rbrace$.
\end{itemize}

This defines a compact Hausdorff topology, thus making $X$ normal. Once again, $\bigcup_{n\in\omega} U_n$ inherits the original topology. Therefore, this latter set is a SHD dense subset of $X$ but $X$ is not SHD due to the countable character of $p$.
\end{exa}
Note that one of the key conditions for a space to be SHD is that it does not have isolated points since every isolated point has a countable local base. Therefore, the first thing we need to verify to determine if a space is SHD is the complete absence of isolated points.\\One of the approaches we considered was to investigate if any compactification or remainder of a compactification of a space is SHD. For example, $\beta\omega$ is not SHD because every point in $\omega$ is isolated in it. However, as we saw previously, $\omega^{*}$ is a SHD space. The key properties of $\omega^{*}$ are that the only convergent sequences in this space are eventually constant, and it has no isolated points. A way to generalize the previous property is as follows using a well known result:
\begin{thm} \label{Coso 1}
Let $X$ be a realcompact non-compact space. Then $X^{*}$ has no isolated points.
\end{thm}
Unfortunately, the previous theorem is not sharp enough in the sense that the converse implication is not true. The following example illustrates this:
\begin{exa}
By \cite[Theorem~8.6.2]{easy} there exists a maximal almost disjoint family $\mathscr{A}$ such that the associated Mrówka space, denoted as $\Psi(\mathscr{A})$, satisfies that $\Psi(\mathscr{A})^{*}$ is homeomorphic to $[0,1]$, which is a compact metric space without isolated points. However, $\Psi(\mathscr{A})$ is not realcompact, as shown in \cite[8H.6]{Gillman}.
\end{exa}
A very useful class of examples is based on the following known result:
\begin{pro} \label{notaalpie}
Let $X$ be a non-compact, locally compact, $\sigma$-compact and Tychnoff space. Then $X^{*}$ is a compact SHD space. 
\end{pro}
This result presents a new way to prove that $\omega^{*}$ is a SHD space. Also, $\left(\mathbb{R}^{n}\right)^{*}$ is a connected SHD space when $n\geq 2$ by \cite[Theorem~1]{conexidad}. This contrasts with the previous examples because they are disconnected SHD spaces. Another result of considerable significance is the following:
\begin{thm}
Let $X$ be a non-compact realcompact space such that the only convergent sequences in $X$ are those that are eventually constant. Then $X^{*}$ is SHD.
\end{thm}
\begin{pf}
Let $p\in X^{*}$. Consider $\{ x_n\}_{n\in\omega}\subseteq X^{*}$ to be an injective sequence converging to $p$. Since $X$ is realcompact, there exists $Z$ which is a zero set in $\beta X$ with $Z\subseteq  X^{*}$ and $Z\cap ({x_n: n\in\omega }\cup\{p\})=\{p\}$. Let $\{ U_n: n\in\omega\}$ be a sequence of cozero sets such that $x_n\in U_n$ and $U_n\cap U_m=\emptyset$ for $n\neq m$, and $p\notin \bigcup_{n\in\omega} U_n$. Hence, $p\in Z\setminus \bigcup_{n\in\omega} U_n=\bigcap_{n\in\omega} (Z_0\setminus U_n)$, and the latter set is a zero set in $\beta X$. Since $\beta X\setminus Z$ is normal (as it is $\sigma$-compact), the function $f:\{x_n: n\in\omega \}\to\mathbb{R}$ defined as
\begin{equation*}
f(x)=\left\{
\begin{array}{cc}
0, & \text{if} \ n \  \text{is even}\\
1, & \text{if} \ n \  \text{is odd}\\
\end{array}
\right.
\end{equation*}
is continuous since $\{x_n: n\in\omega \}$ is closed in $ X^{*}$. Then there exists $F:\beta X\to \mathbb{R}$ that extends $f$. This is impossible as $F$ is continuous and $(x_n)$ converges to $p$. Hence, the only convergent sequences in $X^{*}$ are those that are eventually constant, and since $X$ is realcompact, $ X^{*}$ has no isolated points. It follows that $ X^{*}$ is SHD.
\end{pf}

In the search for more examples of SHD spaces, we came up with the $G_\delta$ modification of a topological space $X$, which involves considering the usual topology of $X$ along with the $G_\delta$ sets as a base for a new topology. We denote the space $X$ with this topology as $X_\delta$. For more information on this topic, see \cite{arangelskiasdasd}. The next theorem sumarizes the work over the $G_\delta$ modification:
\begin{thm} \label{gdelta}
Let $X$ be a regular topological space such that for each $x\in X$, we have $\psi(x,X)>\omega$. If $X_\delta$ denotes the $G_\delta$ modification of $X$, then $X_\delta$ is a SHD, zero-dimensional, and Tychonoff space. Furthermore, if $X$ is also homogeneous, then $X_\delta$ is also homogeneous.
\end{thm}
\begin{pf}
A routine argument shows that $X_\delta$ is a $P$-space without isolated points. Moreover, since $X$ is a regular $P$-space, $X$ is zero-dimensional (see \cite[1W(1)]{Porter}), and thus Tychonoff. By \cite[6L(3)]{Porter}, we have that $X_\delta$ is an $F$-space without isolated points, and therefore, $X_\delta$ is a SHD space.\\\\Finally, if $x,y\in X$ satisfy $x\neq y$, then the homogeneity of $X$ produces a homeomorphism $f : X \to X$ such that $f(x)=y$. It follows that $f : X_\delta \to X_\delta$ is a homeomorphism with $f(x)=y$. Thus, $X_\delta$ is homogeneous. 
\end{pf}
From this theorem, we obtain some interesting examples:
\begin{exa}
Let $\mathfrak{m} > \omega$. It is well known that if $X=D(2)^{\mathfrak{m}}$, then for every $x\in X$, we have $\psi(x,X)>\omega$. By Theorem~\ref{gdelta}, it follows that $X_\delta$ is a SHD, Tychonoff, zero-dimensional, and homogeneous space. But we can say more. Note that $X$ is a topological group, and it can be easily shown using \cite[1W.(7)]{Porter} that $X_\delta$ is also a topological group. Also, if $X_\delta$ were pseudocompact, then by \cite[4AG(6)(e)]{Porter}, $X_\delta$ would be finite, which is absurd. Therefore, $X_\delta$ is a SHD space that is not pseudocompact.
\end{exa}
If we examine the construction in Theorem~\ref{gdelta} in detail, we would expect that if we change the condition ``for every $x\in X$ it holds that $\psi(x,X)>\omega$'' to ``for every $x\in X$ it holds that $\chi(x,X)>\omega$'', the conclusion would still hold. However, the following example shows that this is not possible:
\begin{exa}
By \cite[Chapter 1~11.4]{handbook}, there exists a countable regular space $X$ such that for every $x\in X$ we have $\chi(x,X)=2^{\omega}$. Since $X$ is countable and regular, it follows that $\psi(x,X)=\omega$ for each $x\in X$. As a result, $X_\delta$ is a discrete space, and therefore not SHD.
\end{exa}
A classical and very interesting space in general topology is the Pixley-Roy hyperspace (see \cite{lutzer} and \cite{cardinalpixley}) associated to a $T_1$ topological space $X$, which is the set $\mathscr{F}[X]=\{A\subseteq X : 0<|A|<\omega\}$, i.e., the set of finite subsets of $X$, equipped with the topology generated by local bases for $F\in\mathscr{F}[X]$ of the form $\{[F,V] : V\in\tau_X \text{ and } F\subseteq V\}$, where $[F,V]=\{H\in\mathscr{F}[X] : F\subseteq H\subseteq V\}$. Having presented it, let us observe the following propositions that tell us about the relationship between the SHD property in $X$ and in $\mathscr{F}[X]$. The proof of the next lemma is straightforward: 
\begin{lem} \label{Pixleypespacio}
Let $X$ be a $T_1$ space in which every countable set is closed. If $(x_n)$ is a sequence with no convergent subsequences and $x\in X$, then there exists an open set $U$ in $X$ such that $x\in U$ and moreover, $|U\cap \{x_n: n\in\omega\}|<\omega$.
\end{lem}
Notice that the previous lemma also works for $P$-spaces. 
\begin{pro} \label{Pixleyss1}
Let $X$ be a $T_1$ space where countable sets are closed. If $X$ is SHD, then $\mathscr{F}[X]$ is also SHD.
\end{pro}
\begin{pf}
Let $\{W_n : n \in \omega\}$ be a countable collection of non-empty open sets in $\mathscr{F}[X]$. Since they are non-empty, for each $n \in \omega$ we can choose $F_n \in W_n$. As $W_n$ is open, for each $n \in \omega$ there exists a non-empty open set $U_n$ in $X$ such that $F_n \in [F_n,U_n] \subseteq W_n$. Note that $\{U_n : n \in \omega\}$ is a countable collection of open sets in $X$, and since $X$ is SHD, for each $n \in \omega$ there exists $x_n \in U_n$ such that $(x_n)$ is a sequence without convergent subsequences. For each $n \in \omega$, let's define $H_n = F_n \cup \{x_n\}$. Note that $H_n \in \mathscr{F}[X]$. Moreover, since $x_n \in U_n$ and recalling that $[F_n,U_n] = \{J \in \mathscr{F}[X] : F_n \subseteq H \subseteq U_n\} \subseteq W_n$, then it follows that $F_n \subseteq U_n$, and thus $F_n \subseteq F_n \cup \{x_n\} \subseteq U_n$, i.e., $H_n \in [F_n,U_n] \subseteq W_n$. So, for all $n \in \omega$, it holds that $H_n \in W_n$.\\\\Let's prove that the sequence $(H_n)\subseteq \mathscr{F}[X]$ has no convergent subsequences. For this, let's consider $(H_{n_k})$ as a subsequence of $(H_n)$ along with a fixed $H\in \mathscr{F}[X]$. Note that the subsequence $(H_{n_k})$ induces a subsequence $(x_{n_k})$ of $(x_n)$. Since this sequence has no convergent subsequences, then $(x_{n_k})$ also has no convergent subsequences. Since $H$ is finite, we can label its elements without repetitions so that $H=\{z_1,\dots,z_\ell\}$. Thus, for each $i\in\{1,\dots,\ell\}$, by Lemma~\ref{Pixleypespacio}, there exists an open set $O_i$ in $X$ such that $z_i\in O_i$ and $|O_i\cap\{x_{n_k}: k\in\omega\}|<\omega$. If we consider $O=\bigcup\{O_i: i\in\{1,\dots,\ell\}\}$, then $O$ is an open set such that $H\subseteq O$. Let's consider the neighborhood $[H,O]$. Thus, note that no tail of the sequence $(H_{n_k})$ is contained in $[H,O]$ since $O$ only contains a finite number of terms from $(x_{n_k})$. Thus, $(H_{n_k})$ does not converge to $H$. Therefore, $\mathscr{F}[X]$ is a SHD space.
\end{pf}
\begin{cor} \label{ppixley}
Let $X$ be a $T_1$ space and a $P$-space. If $X$ is SHD, then $\mathscr{F}[X]$ is SHD.
\end{cor}
The following result can be proven with standard arguments:
\begin{pro} \label{equivalenciapespaciopixley}
Let $X$ be a $T_1$ space. Then $X$ is a $P$-space if and only if the Pixley-Roy hyperspace $\mathscr{F}[X]$ is also a $P$-space.
\end{pro}
\begin{thm}
Let $X$ be a regular topological space. Then, the Pixley-Roy hyperspace $\mathscr{F}[X]$ is a $P$-space and SHD if and only if $X$ is a $P$-space and SHD.
\end{thm}
\begin{pf}
If $X$ is a $P$-space and SHD, by Corollary~\ref{ppixley} and Proposition~\ref{equivalenciapespaciopixley}, it follows that $\mathscr{F}[X]$ is a $P$-space and SHD. Furthermore, if $\mathscr{F}[X]$ is a $P$-space, then $X$ is a $P$-space by Proposition~\ref{equivalenciapespaciopixley}. Since $\mathscr{F}[X]$ is SHD, it does not have isolated points. Hence, $X$ does not have isolated points either. Thus, $X$ is a $P$-space without isolated points, and being regular, it is an $F$-space without isolated points by \cite[6L]{Porter} and therefore $X$ is SHD.
\end{pf}

\section{SHD spaces of all cardinalities}

The purpose of this section is to show that there exist SHD spaces of all infinite cardinalities with various topological properties. We will be using the following result constantly (see Proposition~\ref{Fespacios}).

\begin{thm}\label{thm_ed_SS1} If $X$ is an extremally disconnected Tychonoff space with no isolated points, then $X$ is SHD.

\end{thm}

Let $X$ be a Hausdorff space. The symbol $\theta X$ will stand for the collection of all $\R(X)$-ultrafilters in $X$. Furthermore, for each $A\in \R(X)$ let us denote by $\lambda(A) := \{\mathscr{U} \in \theta X : A\in\mathscr{U}\}$. It can be verified that the family $\{\lambda(A) : A\in \R(X)\}$ is a basis for a topology on $\theta X$.

Let us use the symbol $EX$ to refer to the subspace $\{\mathscr{U} \in \theta X : \bigcap\mathscr{U} \neq \emptyset\}$ (i.e., the {\it Iliadis absolute} of $X$) and, for each $x\in X$, let $F(x):=\{A\in \R(X) : x\in \int_X A\}$. It is known that the function $k_X : EX \to X$ determined by $k_X(\mathscr{U}) \in \bigcap\mathscr{U}$ is well defined and has the following properties which can be found in \cite[Theorem~(e), p.~459]{Porter}.

\begin{thm}\label{thm_EX} Let $X$ be a Hausdorff space.

\begin{enumerate}
\item $EX$ is an extremally disconnected zero-dimensional Hausdorff space.
\item If $\mathscr{U} \in \theta X$ and $x\in X$, then $\mathscr{U} \in EX$ and $k_X(\mathscr{U})=x$ if and only if $F(x)\subseteq \mathscr{U}$.
\end{enumerate}

\end{thm}

All that remains is to mention the following lemma that appears in \cite[Theorem~(b), p.~445]{Porter}.

\begin{lem}\label{lemma_dense_ed} If $D$ is a dense subspace of an extremally disconnected space, then $D$ is extremally disconnected.

\end{lem}

With this background we can prove the following result.

\begin{thm}\label{thm_1} If $X$ is a Hausdorff space with no isolated points, then there exists an space $sX$ that satisfies the following conditions:

\begin{enumerate}
\item $sX$ is an extremally disconnected zero-dimensional Hausdorff space without isolated points (consequently, SHD).
\item $|X|=|sX|$, $\pi w(X) = \pi w(sX)$ and $c(X) = c(sX)$.
\end{enumerate}

\end{thm}

\begin{pf} For each $x\in X$ fix $\mathscr{U}_x \in k_X^{-1}\{x\}$ and consider $sX := \{ \mathscr{U}_x : x\in X\}$ with the topology that it inherits as a subspace of $EX$. Clearly, $k_X$ is a bijection between $sX$ and $X$.

\medskip

\noindent {\bf Claim.} $sX$ is a dense subspace of $EX$ with no isolated points.

\medskip

To verify the previous statement, it is enough to prove that if $A\in \R(X)\setminus\{\emptyset\}$, then $|\lambda(A)\cap sX|\geq 2$. Let $A$ be a regular non-empty closed subset of $X$. First, since $A$ is regular closed, $\int_X A$ is not empty and therefore, since $X$ has no isolated points, distinct $x,y\in \int_X A$ exist. Then, since $k_X(\mathscr{U}_x) = x$ and $k_X(\mathscr{U}_y)=y$, Theorem~\ref{thm_EX}(2) implies that $F(x)\subseteq \mathscr{U}_x$ and $F(y)\subseteq\mathscr{U}_y$. Thus, since $A \in F(x)\cap F(y)$, we deduce that $A\in \mathscr{U}_x\cap \mathscr{U}_y$ and hence $\{ \mathscr{U}_x, \mathscr{U}_y\} \subseteq \lambda(A)$; in particular, $|\lambda(A)\cap sX|\geq 2$.

Now, since $sX$ is a dense subspace of $EX$, we obtain that $sX$ is Hausdorff, zero-dimensional and extremally disconnected (see Theorem~\ref{thm_EX}(1) and Lemma~\ref{lemma_dense_ed}). Thus, Theorem~\ref{thm_ed_SS1} and our Claim guarantee that $sX$ is SHD.

Finally, to verify that $\pi w(X) = \pi w(sX)$ and $c(X) = c(sX)$, we only have to remember a couple of results. Recall that if $D$ is a dense subspace of a $T_3$ space $Y$, then $\pi w(D) = \pi w(Y)$ and $c(D) = c(Y)$ (see \cite[2.6(a), p.~14]{juhasz1980} and \cite[2.7(a), p.~15]{juhasz1980}). Furthermore, for any Hausdorff space $Y$ it is satisfied that $\pi w(Y) = \pi w(EY)$ and $c(Y) = c(EY)$ (see \cite[6B(4), p.~496]{Porter}). In short, since $X$ is $T_2$, $sX$ is a dense subspace of $EX$ and $EX$ is $T_3$, the relations $$\pi w(X) = \pi w(EX) = \pi w(sX) \quad \text{y} \quad c(X) = c(EX) = c(sX)$$ are verified.
\end{pf}

Since for any infinite cardinal $\kappa$ it is satisfied that the free topological sum $\bigoplus_{\alpha<\kappa} \mathbb{Q}$ is a Hausdorff space with no isolated points having cardinality, $\pi $-weight and cellularity $\kappa$, Theorem~\ref{thm_1} implies:

\begin{cor} For any infinite cardinal $\kappa$ there exists a Hausdorff space $X_\kappa$ that is zero-dimensional, SHD and has cardinality, $\pi$-weight and cellularity $\kappa$.

\end{cor}

In particular, a fundamental example that we will constantly mention and that is obtained from what is stated in Theorem~\ref{thm_1} is the following:

\begin{exa}\label{ej_sQ} The space $s\mathbb{Q}$ is $T_2$, zero-dimensional, extremally disconnected (consequently, SHD), countable, and has countable $\pi$-weight.

\end{exa}

What follows is to expose the details to obtain a result similar to Theorem~\ref{thm_1} related to the cardinal function known as density.

If $X$ is a Hausdorff space, the space $\theta X$ is compact, $T_2$, extremally disconnected, zero-dimensional, and if $X$ has no isolated points, then $\theta X$ also has no isolated points (see \cite[\S 6.3]{Porter}). Furthermore, if $X$ is infinite a consequence of \cite[Theorem~B.14, p.~270]{comneg1982} is that $w(\theta X)=|\R(X)|$.

In the case of compact Hausdorff spaces, the density of the space $\theta X$ coincides with the density of the original space. This fact is proven in \cite[Corollary~B.20, p.~272]{comneg1982}.

The previous remarks combined with Theorem~\ref{thm_ed_SS1} imply the following result.

\begin{thm}\label{thm_2} If $X$ is an infinite Hausdorff space with no isolated points, then $\theta X$ is a zero-dimensional SHD compact Hausdorff space with $w(\theta X) = |\R(X) |$. Also, if $X$ is compact, $d(X)=d(\theta X)$.

\end{thm}

With the above theorem we are positioned to obtain a result similar to Theorem~\ref{thm_1} with respect to density.

\begin{thm}\label{thm_3} For any infinite cardinal $\kappa$, there exists a space $Y_\kappa$ that is Hausdorff, compact, zero-dimensional and SHD with $d(Y_\kappa)=2^{\kappa} \leq w( Y_\kappa) \leq 2^{2^{\kappa}}$.

\end{thm}

\begin{pf} Let $\lambda$ stand for the cardinal $2^{\kappa}$. We shall use the symbol $\lambda^{*}$ to represent the residue of the Stone-\v{C}ech compactification of the discrete space of cardinality $\lambda$. By Theorem~\ref{thm_2} it is only necessary to argue that the space $Y_\kappa := \theta\left( \lambda^{*}\right)$ satisfies the relations $d(Y_\kappa)=\lambda \leq w( Y_\kappa) \leq 2^{\lambda}$.

On the one hand, \cite[Theorem, p.~229]{walker1974} implies that $d(Y_\kappa)=d\left(\lambda^{*}\right)=\lambda^{\omega}=\lambda$. On the other hand, by virtue of the equality $w(Y_\kappa)=\left|\R(\lambda^{*})\right|$, it is enough to notice that usual arguments with cardinal functions (see \cite{handbook}) show that $\lambda=d\left(\lambda^{*}\right)\leq w\left(\lambda^{*}\right)\leq \left|\R\left(\lambda^{ *}\right)\right|\leq 2^{d\left(\lambda^{*}\right)}=2^{\lambda}$.
\end{pf}

For any infinite cardinal $\kappa$, let us denote by $D(2)^{\kappa}$ the Cantor cube of weight $\kappa$. A routine argument shows that $D(2)^{\kappa}$ is a compact Hausdorff space with no isolated points. Furthermore, a consequence of the Generalized Continuum Hypothesis, \textsf{GCH}, and of \cite[Example~11.8, p.~44]{handbook} is that \[d\left(D(2)^{2^ {\kappa}}\right)=\kappa \quad \text{y} \quad \left|\R\left(D(2)^{2^{\kappa}}\right)\right|=2^ {\kappa}.\]

\begin{thm}\label{thm_4} $[\mathsf{GCH}]$ For any infinite cardinal $\kappa$, there exists a space $Z_\kappa$ that is Hausdorff compact, zero-dimensional and SHD with $d(Z_\kappa)=\kappa$ and $w(Z_\kappa )=2^{\kappa}$.

\end{thm}

\begin{pf} By virtue of Theorem~\ref{thm_2} and the observations in the previous paragraph, if $\kappa\geq\omega$, it is clear that $Z_\kappa := \theta\left( D(2)^{2^{\kappa }}\right)$ satisfies the desired properties.
\end{pf}

\section{Non-semiregular SHD spaces}

Although topological constructions are usually made in such a way that the resulting space possesses rich properties, we will show in this section that SHD spaces are also very versatile in the sense that SHD spaces can be found, again of all the possible infinite cardinalities, and satisfying not being semiregular (i.e., the set of regular open sets is not a basis for the corresponding spaces).

Recall that a space $X$ admits {\bf no} non-trivial convergent sequences if and only if the only convergent sequences in $X$ are semiconstant. Let's recall Proposition~\ref{lemma_SS1}, which has the following Corollary. These results will be used repeatedly in the development of this section.

\begin{cor}\label{cor_SS1} If $X$ is a $T_1$ space without isolated points that does not admit non-trivial convergent sequences, then $X$ is SHD.

\end{cor}

Next, we will present a general way of producing non-semiregular SHD topological spaces from topological spaces with certain characteristics.

\begin{definition} If $X$ is a topological space, we will denote by $mX$ the set $X$ equipped with the topology whose base is the collection $$\mathscr{B}:= \left\{U\setminus A : U\in \tau_X \ \wedge \ A\in [X]^{\leq\omega}\right\}.$$

\end{definition}

\begin{thm}\label{thm_5} If $X$ is a topological space such that, for any $U\in\tau_X^+$, it is satisfied that $|U| > \omega$, then $mX$ is SHD. Furthermore, if $X$ is $T_0$, $T_1$, $T_2$, Urysohn, or completely Hausdorff, then so is $mX$.

\end{thm}

\begin{pf} First note that if $U\in \tau_{mX}^{+}$, then there exist $V\in \tau_X^+$ and $A\in [X]^{\leq\omega}$ with $V\setminus A \subseteq U$. Then, since $|V|>\omega$ we deduce that $|U|\geq |V\setminus A| >\omega$; in particular, $U$ is infinite.

To prove that $mX$ does not admit non-trivial convergent sequences, we note that if $x\in mX$ and $(x_n)$ is a sequence in $mX$ with $x\not \in \{x_n : n<\omega \}$, then $U:= X\setminus \{x_n : n<\omega\}$ is an element of $\tau_{mX}(x)$ such that $U\cap \{x_n : n <\omega\} = \emptyset$; consequently, $(x_n)$ does not converge to $x$ in $mX$ and, therefore, any convergent sequence in $mX$ is necessarily semiconstant. Consequently, Proposition~\ref{lemma_SS1} guarantees that $mX$ is SHD.

For the second part it is enough to observe that if $X$ has any of the separation properties mentioned in the statement of our theorem, then the inclusion $\tau_X \subseteq \tau_{mX}$ guarantees that $mX$ satisfies the same property.
\end{pf}

\begin{lem}\label{lemma_mX} If $X$ is a topological space and for any $U\in\tau_X^+$ we have that $|U| > \omega$, then $\RO(mX) \subseteq \tau_X$. In particular, if there exists $V\in \tau_{mX} \setminus \tau_X$, then $mX$ is not semiregular.

\end{lem}

\begin{pf} Let $U\in \RO(mX)$, $x\in U$, $V\in \tau_X$ and $ A\in [X]^{\leq\omega}$ be such that $x\in V\setminus A \subseteq U$. We will prove first that $\cl_{mX} (V\setminus A) = \cl_X V$.

Start by noticing that the relations $\tau_X \subseteq \tau_{mX}$ and $V\setminus A \subseteq V$ imply that $\cl_{mX} (V\setminus A) \subseteq \cl_{mX} V \subseteq \cl_X V$. On the other hand, if $y\in \cl_X V$, $W\in \tau_X$ and $B\in [X]^{\leq\omega}$ are such that $y\in W\setminus B$, then $W\cap V\in \tau_X^{+}$. Thus, it is satisfied that $|W\cap V|>\omega$ and, therefore, we obtain that $V \cap (W\setminus B) \neq \emptyset$. Hence, $y\in \cl_{mX} V$.

To finish the argument we note that $V\in \tau_X(x)$ (in particular, $V\in \tau_{mX}(x)$) and $V\subseteq \cl_X V = \cl_{mX} (V \setminus A)$. Thus, $x\in V \subseteq \int_{mX} \cl_{mX} (V\setminus A) \subseteq \int_{mX} \cl_{mX} U =U$; in short, $x\in V \subseteq U$.

The second part is simple: if $mX$ is semiregular, then any open set $mX$ is a union of elements of $\RO(mX)$ which, by the inclusion we proved above, turns out to be a union of elements of $\tau_X$, i.e., it is open in $X$.
\end{pf}

\begin{lem}\label{lemma_no_num} There is a subspace $X$ of $\mathbb{R}$ with the following characteristics:
\begin{enumerate}
\item $|X| = \omega_1$;
\item for any $U\in\tau_{X}^{+}$, $|U|=\omega_1$ and $U\cap \mathbb{Q}\neq \emptyset$; and
\item $X\setminus \mathbb{Q}$ is not an open subset of $X$.
\end{enumerate}

\end{lem}

\begin{pf} Let $\mathscr{B}:=\{(a,b) : a,b\in\mathbb{Q} \ \wedge \ a<b\}$. For each $B\in\mathscr{B}$ choose $X_B\in [B]^{\omega_1}$ with $X_B\cap \mathbb{Q} \neq \emptyset$, and consider the subspace $X:=\bigcup\{X_B : B\in\mathscr{ B}\}$. Clearly, $|X| = \omega_1$. Furthermore, if $U\in\tau_{\mathbb{R}}^{+}$ and $B\in \mathscr{B}$ are such that $B\subseteq U$, then we have the relations \begin{align*} \omega_1 = |X_B \cap B| \leq |X\cap B| \leq |X\cap U| \leq |X|\leq |\mathscr{B}|\cdot \sup\left\{|X_B| : B\in\mathscr{B}\right\} \leq\omega_1.
\end{align*} In conclusion, $|X\cap U|=\omega_1$ and $(X\cap U)\cap \mathbb{Q} \neq \emptyset$. Finally, since (1) and (2) ensure that $X\cap \mathbb{Q}$ is a proper and dense subset of $X$, we infer that $X\setminus \mathbb{Q}$ does not belong to the collection $\tau_X$.
\end{pf}

\begin{thm}\label{thm_6} For any uncountable cardinal $\kappa$ there exists a space $X_\kappa$ that is completely Hausdorff, not semiregular, and SHD of cardinality $\kappa$.

\end{thm}

\begin{pf} Let $X$ be as in Lemma~\ref{lemma_no_num}. For each $\alpha<\kappa$ let us denote by $Y_\alpha$ the space $X\times\{\alpha\}$. We will use the symbol $Y$ to refer to the topological sum $\bigoplus_{\alpha<\kappa} Y_\alpha$. Finally, let $X_\kappa := mY$. Observe that since Lemma~\ref{lemma_no_num}(1) ensures the equality $|X|=\omega_1$, it is verified that $X_\kappa$ has cardinality $\kappa$.

Now, since $X$ is a subspace of $\mathbb{R}$, $\mathbb{R}$ is completely Hausdorff, and this property is hereditary and is preserved under topological additions, we get that $Y$ is a completely Hausdorff space. On the other hand, Lemma~\ref{lemma_no_num}(2) guarantees that $|U|>\omega$ provided that $U\in \tau_Y^+$. Thus, Theorem~\ref{thm_5} certifies that $X_\kappa$ is completely Hausdorff and SHD.

To verify that $X_\kappa$ is not semiregular let us first observe that, by Lemma~\ref{lemma_no_num}(3), $(X\setminus \mathbb{Q})\times \{0\}$ is not a open subset of $Y_0$; in particular, $(X\setminus \mathbb{Q})\times \{0\}$ is not an element of $\tau_Y$. However, since \[(X\setminus \mathbb{Q})\times \{0\}= (X\times \{0\}) \setminus (\mathbb{Q}\times \{0\}) ,\] the membership $(X\setminus \mathbb{Q})\times \{0\}\in \tau_{mY}$ is satisfied. Thus, Lemma~\ref{lemma_mX} implies that $X_\kappa$ is not semiregular.
\end{pf}

The rest of the section is dedicated to constructing a countable space which is SHD and not semiregular with a method completely different from the construction exposed in Theorem~\ref{thm_6}.

Recall that a Hausdorff space $X$ is {\it $H$-closed} if it is closed in any Hausdorff space that contains it as a subspace. On the other hand, a function between topological spaces $f:X\to Y$ is {\it $\theta$-continuous} if for any $x\in X$ and $V\in \tau_Y(f(x))$, there exists $U\in \tau_X(x)$ such that $f[\cl_X U]\subseteq \cl_Y V$. Our next lemma requires familiarity with the construction of Theorem~\ref{thm_1}.

\begin{lem}\label{lemma_not_H-closed} Let $X$ be a topological space. If $A$ is a clopen subspace of $X$ that is not $H$-closed, then $\lambda(A)\cap sX$ is not $H$-closed.

\end{lem}

\begin{pf} We will first argue that $k_X[\lambda(A)\cap sX] = A$. Note that \cite[Theorem~(e)(3), p.~459]{Porter} guarantees the relations $k_X[\lambda(A)\cap sX] \subseteq k_X[\lambda(A)\cap EX] =A$. On the other hand, if $x\in A$, then $x\in \int_X A$ and therefore $A\in F(x)$. Then, since $k_X(\mathscr{U}_x) = x$, Theorem~\ref{thm_EX}(2) guarantees that $F(x) \subseteq \mathscr{U}_x$. In this way, $A \in \mathscr{U}_x$, that is, $\mathscr{U}_x \in \lambda(A) \cap sX$ and thus, $x\in k_X[\lambda(A )\cap sX]$. In conclusion, $k_X[\lambda(A)\cap sX] = A$.

Finally, since $k_X$ is a $\theta$-continuous function (see \cite[Theorem~(e)(5), p.~459]{Porter}) and $A$ is not $H$-closed, \cite[Proposition~(h), p.~302]{Porter} implies that $\lambda(A)\cap sX$ is not $H$-closed.
\end{pf}

If $X$ is a Hausdorff space and $\mathscr{F}$ is an open filter base on $X$, then $\mathscr{F}$ is free if and only if the {\it adherence of $\mathscr{F}$}, $a_X(\mathscr{F}) := \bigcap\{\cl_X F : F\in\mathscr{F}\}$, is empty (see \cite[\S2.3]{Porter}). Additionally, we will use the symbol $\kappa X$ to denote the {\it Kat\v{e}tov extension of $X$}, that is, $\kappa X$ is the set $\kappa X := X \cup \{\mathscr{U} : \mathscr{U} \ \text{is a free open ultrafilter on} \ X\}$ whose topology is determined by the basis $\tau_X \cup \{U \cup \{ \mathscr{U}\} : \mathscr{U}\in \kappa X\setminus X \ \text{y} \ U\in \mathscr{U}\}$ (see \cite[\S4.8]{Porter}).

\begin{lem}\label{lemma_ultrafilter} If $X$ is a Hausdorff space and $U\in \tau_X$, then $\cl_X U$ is not $H$-closed if and only if there exists $\mathscr{U}\in \kappa X\setminus X$ with $U\in \mathscr{U}$.

\end{lem}

\begin{pf} Let $Y$ denote the space $\cl_X U$ and suppose first that there exists $\mathscr{U}\in \kappa X\setminus X$ with $U\in \mathscr{U}$. Consider the collection $\mathscr{F} := \{U\cap V : V\in\mathscr{U}\}$. A routine argument shows that $\mathscr{F}$ is an open filter base on $Y$. Also, since $a_Y(\mathscr{F}) \subseteq \bigcap\{\cl_X (U\cap V) : V\in\mathscr{U}\} \subseteq a_X(\mathscr{U})=\emptyset $, we deduce that $\mathscr{F}$ is free, and therefore any open filter on $Y$ that extends $\mathscr{F}$ must be free. By \cite[Proposition~(b), p.~298]{Porter}, we deduce that $Y$ is not $H$-closed.

Suppose now that $Y$ is not $H$-closed and use \cite[Proposition (b), p.~298]{Porter} to find an ultrafilter $\mathscr{W}$ that is open and free in $Y $. Note that $U\in\mathscr{W}$ because, otherwise, since $\mathscr{W}$ is an open ultrafilter in $Y$, we would have $\emptyset = Y \setminus \cl_Y U \in \mathscr{W}$, which is absurd. Furthermore, since $Y$ is a closed subset of $X$, for any $W\in \mathscr{W}$ the equality $\cl_Y W = \cl_X W$ is satisfied (see \cite[1A(2), p.~55]{Porter}).

Now, consider the collection $\mathscr{F}:=\{V\in \tau _X : V\cap Y \in \mathscr{W}\}$. A simple reasoning shows $\mathscr{F}$ is an open filter on $X$ that satisfies the membership $U\in\mathscr{F}$. For this reason, there exists an open ultrafilter $\mathscr{U}$ in $X$ such that $\mathscr{F} \subseteq \mathscr{U}$ (see \cite[Proposition~(d)(2), p.~93]{Porter}). In order to see that $a_X(\mathscr{U})=\emptyset$ suppose, in search of a contradiction, that $x\in a_X(\mathscr{U})$. Given $W\in\mathscr{W}$ there exists $V\in\tau_X$ such that $W=V\cap Y$. In particular, $V\in \mathscr{F} \subseteq \mathscr{U}$. Then, as $U\in \mathscr{U}$ we have $V\cap U\in \mathscr{U}$ and therefore, $x\in \cl_X(V\cap U)$. Finally, we can use \cite[Proposition (a)(3), p.~81]{Porter} to ensure the equalities \[\cl_X(V\cap U) = \cl_X(V\cap \cl_X U)= \cl_X(V\cap Y)= \cl_X W = \cl_Y W.\] This shows that for any $W\in\mathscr{W}$, $x\in \cl_Y W$; in other words, $x\in a_Y(\mathscr{W})$, which contradicts that $\mathscr{W}$ is free. In sum, $\mathscr{U}\in \kappa X\setminus X$ and $U\in \mathscr{U}$.
\end{pf}

In what follows we will be working with the space $s\mathbb{Q}$ (see Example~\ref{ej_sQ}).

\begin{lem}\label{lemma_ultrafilters_sQ} There exists a collection $\{\mathscr{V}_n : n<\omega\}\subseteq \kappa\left(s\mathbb{Q}\right)\setminus s\mathbb{Q}$ such that, for for any $U\in \tau_{s\mathbb{Q}}^+$, there exists $n<\omega$ with $U\in \mathscr{V}_n$.

\end{lem}

\begin{pf} Since $\{(a,b)\cap \mathbb{Q} : a,b\in \mathbb{R}\setminus\mathbb{Q} \ \wedge \ a<b\}$ is a basis for $\mathbb{Q}$ and this space is second countable, the family admits a subset $\{B_n : n<\omega\}$ enumerated faithfully that is a basis for $\mathbb{Q}$. Note that if $n<\omega$ then $B_n$ is a clopen subset of $\mathbb{Q}$. Furthermore, since $B_n$ is not closed in $\mathbb{R}$, we deduce that $B_n$ is not $H$-closed. For this reason, for each $n<\omega$ we can use Lemma~\ref{lemma_not_H-closed} to verify that $\lambda(B_n)\cap s\mathbb{Q}$ is not $H$-closed. Thus, since $\lambda(B_n)\cap s\mathbb{Q}$ is a clopen subspace of $s\mathbb{Q}$ and $\cl_{s\mathbb{Q}} \left(\lambda (B_n)\cap s\mathbb{Q}\right)$ is not $H$-closed, Lemma~\ref{lemma_ultrafilter} implies the existence of $\mathscr{V}_n \in \kappa\left(s \mathbb{Q}\right)\setminus s\mathbb{Q}$ with $\lambda(B_n)\cap s\mathbb{Q} \in \mathscr{V}_n$.

Finally, if $U\in \tau_{s\mathbb{Q}}^{+}$, then there exists $A\in \R(\mathbb{Q})\setminus\{\emptyset\}$ such that $\lambda(A)\cap s\mathbb{Q}\subseteq U$. Then, since $\int_\mathbb{Q} A$ is non-empty, there exists $n<\omega$ with $B_n \subseteq \int_\mathbb{Q} A$. Thus, the inclusion $\lambda(B_n)\cap s\mathbb{Q}\subseteq \lambda(A)\cap s\mathbb{Q}$ certifies that $U\in \mathscr{V}_n$. 
\end{pf}

In the Example~\ref{ej_sQ} it is mentioned that the space $s\mathbb{Q}$ has the Hausdorff property, is extremally disconnected and has no isolated points. We will use these facts in the proof of Theorem~\ref{thm_7}. For the purposes of the following result, the symbol $\mathbb{Q}_{*}$ will represent the subspace $s\mathbb{Q} \cup \{\mathscr{V}_n : n <\omega\}$ of $\kappa\left(s\mathbb{Q}\right)$, where $\{\mathscr{V}_n : n <\omega\}$ is as in Lemma~\ref{lemma_ultrafilters_sQ}.

\begin{thm}\label{thm_7} $\mathbb{Q}_{*}$ is Hausdorff, extremally disconnected, without isolated points and not semiregular. In particular, $\mathbb{Q}_{*}$ is a countable space, SHD and not semiregular.

\end{thm}

\begin{pf} First, since $s\mathbb{Q}$ has no isolated points and is dense in $\mathbb{Q}_{*}$, $\mathbb{Q}_{*}$ also does not have isolated points. On the other hand, since $s\mathbb{Q}$ is extremally disconnected, \cite[Theorem~(b)(7), p.~445]{Porter} guarantees that $\kappa\left(s\mathbb{Q }\right)$ is too. Thus, since $\mathbb{Q}_{*}$ is dense in $\kappa\left(s\mathbb{Q}\right)$, \cite[Theorem~(b)(2), p. ~445]{Porter} ensures that $\mathbb{Q}_{*}$ is extremally disconnected and has the Hausdorff property.

To argue that $\mathbb{Q}_{*}$ is not semiregular let's remember that, in extremally disconnected spaces, semiregularity and zero-dimensionality coincide (see \cite[Theorem, p.~451]{Porter}). With this idea in mind, it is enough to verify that $\mathbb{Q}_{*}$ is not completely regular.

Suppose, in search of an absurdity, that $\mathbb{Q}_{*}$ is a completely regular space. Let us fix $x\in s\mathbb{Q}$ and use that $\{\mathscr{V}_n : n<\omega\}$ is a closed subspace of $\mathbb{Q}_{*}$ (see \cite[Proposition~(p)(1), p.~309]{Porter}) to find a continuous function $f: \mathbb{Q}_{*} \to [0,1]$ such that $f(x)=0$ and $f\left[\{\mathscr{V}_n : n<\omega\}\right] = \{1\}$. Let us use the continuity of $f$ to find $U\in \tau_{s\mathbb{Q}}(x)$ with $f[U] \subseteq [0,1/2)$. Hence, since $\mathbb{Q}_{*}$ is regular and $U\in \tau_{\mathbb{Q}_{*}}(x)$, there exists $V\in \tau_{\mathbb{ Q}_{*}}(x)$ with $\cl_{\mathbb{Q}_{*}} V \subseteq U$.

Now, note that the inclusion $V\subseteq U$ guarantees that, in fact, $V\in \tau_{s\mathbb{Q}}^{+}$. Thus, Lemma~\ref{lemma_ultrafilters_sQ} produces $n<\omega$ with $V\in \mathscr{V}_n$. Also, as $\mathscr{V}_n \in \cl_{\kappa\left(s\mathbb{Q}\right)} V$ (see \cite[Proposition~(p)(2), p.~309 ]{Porter}), we deduce that $\mathscr{V}_n \in \cl_{\mathbb{Q}_{*}} V$. Finally, since $f\left[\cl_{\mathbb{Q}_{*}} V\right] \subseteq [0,1/2)$ and $f(\mathscr{V}_n)=1$, we have obtained the desired contradiction. Consequently, $\mathbb{Q}_{*}$ is not semiregular.
\end{pf}

\section{Open questions}
\begin{question} Is it possible to obtain a result similar to Theorem~\ref{thm_4} without additional axioms?

\end{question}
\begin{question}
Is it true that if $\kappa$ is an uncountable cardinal, then $X=D(2)^{\kappa}$ is a SHD space?
\end{question}
In the same way of Theorem~\ref{Coso 1} we have the next question: 
\begin{question}
If $X$ is realcompact, Tychonoff and non-compact, is it true that $ X^{*}$ is SHD?
\end{question}
As we have seen in \ref{Ejemplo3}, the property of being SHD is not preserved under normal extensions; however, the Stone-Čech compactification is a very particular extension with strong properties. In this regard, we have the following questions:
\begin{question}
If $X$ is Tychonoff, non-compact and SHD, does it hold that $\beta X$ is SHD? 
\end{question}
\begin{question}
The property of being SHD is dense hereditary?
\end{question}
\begin{question}
If $X$ is Tychonoff and non-compact, is it always true that if $\beta X$ is SHD then $X$ is SHD? 
\end{question}
\begin{question}
Is $\mathscr{F}[X]$ SHD whenever $X$ is SHD and $T_1$?
\end{question}
\begin{question}
Is it true that if $X$ is $T_1$ and $\mathscr{F}[X]$ is SHD, then $X$ is SHD? 
\end{question}

\end{document}